\documentclass[12pt]{article}
\usepackage{amssymb,amsmath}

\newtheorem{thm}{Theorem}

\newtheorem{lemma}[thm]{Lemma}

\newenvironment{remk}{\noindent{\sc
Remark~:}}{\goodbreak\vskip10pt}

\def\carrebb{\thinspace\hbox{\vrule height 6pt\vbox{
\hrule width 5pt\kern 5pt\hrule width
5pt}\vrule}\thinspace}\def\carreb{\thinspace\hbox{\vrule height 6pt\vbox{
\hrule width 5pt\kern 5pt\hrule width 5pt}\vrule}\thinspace}
\def\build#1_#2^#3{\mathrel{\mathop{\kern 0pt#1}\limits_{#2}^{#3}}}

\def\R{{\bf R}}

\def\Z{{\bf Z}}
\def\S{{\bf S}}
\begin{document}

\thispagestyle{empty}
\title{ $C^1$-generic billiard tables have a dense set of periodic points}
\author{Marie-Claude ARNAUD
\thanks{Avignon University, LMA, EA 2151   F-84000, Avignon, France 
.}
 \thanks{Membre de l'Institut universitaire de France}}
\maketitle

\abstract{We prove that the set of periodic points of a generic
$C^1$-billiard table  is dense in the
phase space.}
\bigskip

\tableofcontents
\newpage

\section{Introduction}

\rm
 In the general domain of conservative dynamical systems, periodic orbits have
an important place. Some questions concerning them are :\medskip

\bf 1.\quad\rm \rm Do they exist ?\medskip

\bf 2.\quad\rm  What is their type (are they elliptic, hyperbolic...)?\medskip

\bf 3.\quad\rm  How many are they (for example, what is the number of
periodic orbits
having  less than $N$ bounces?) ?\medskip

\bf 4.\quad \rm Is the set of periodic points dense in the considered
manifold ?\dots etc \dots\medskip

In the particular domain of smooth strictly convex billiard tables, some
answers
concerning the question 1)  (the proof is due to Birkhoff ; see for example
\cite{KT}), the
question 2) (see \cite{KT}) and the question 3) (see \cite{M}) exist.
Some answer
to 1)   in the case of smooth non convex billiard is given in \cite{BG}. But
nowhere the case
of the density of periodic points in smooth tables is treated. To be complete, let us mention an example of a $C^\infty$ billiard table in \cite{Tr} where the dynamics is formally conjugated to an irrational rigid rotation on some open set: if the conjugacy could be a true conjugacy, the billiard map has an open set without periodic point.\medskip

However, in the domain of conservative dynamical systems (symplectic
diffeomorphisms for
example), it has been proved in \cite{PR} (see \cite{A}
 too for a slightly different
proof) that
every $C^1$-generic symplectic diffeomorphism of a compact manifold has a
dense subset
of periodic points. This result is the consequence of a hard theorem,
called the
``$C^1$-closing lemma''. The $C^0$-version of this result is quite easier.
Let us explain
in   few words the proof of the ``$C^0$-closing lemma ''~:\medskip

\noindent{\bf $C^0$-closing lemma }: {\sl Let $f : M \longrightarrow M$ be a
homeomorphism of a manifold
$M$ and $x\in M$ a positively recurrent point  for $f$. Let ${\cal U}$ be a
$C^0$-
neighbourhood of $f$. Then there exists $g\in {\cal U}$ such that $x$ is
periodic for
$g$.}

\medskip

\noindent {\sl Proof~:} \rm \quad There exists a connected open neighbourhood $V$
of
$x$ such that
:

``{\sl if $g : M\longrightarrow M$ is a homeomorphism and if  $support(g\circ
f^{-1})
\subset V$ then $g\in {\cal U}$}''\medskip

The point $x$  being recurrent, there exists $N \geq 1$ such that $f^N x \in V$. We choose
$N$ as small as possible. There are two
cases :
\begin{enumerate}
\item[--] $x$ is periodic for $f$ ; we choose $g=f$.\medskip
\item[--]$x$ is not periodic for $f$. Then   $V$
is a connected
neighbourhood of $x$ and $f^Nx$. There exists $h : M \longrightarrow M$
homeomorphism
such that $h(f^Nx) = x$ and $support(h)\subset V$.\medskip

\noindent If we define $g$ by $g = h\circ f$, then :\medskip
\begin{enumerate}
\item[$\bullet  $] $support(g\circ f^{-1}) = \; support(h) \subset V$, thus $g\in
{\cal U}$
;
\item[$\bullet$] $f(x) \notin \; support(h)$,\dots, $f^{N-1}(x) \notin\;
support(h)$. Thus $g(x)
= h\circ f(x) = f(x)$, \dots , $ g^{N-1}(x) = g\circ g^{N-2}(x) =
g\big(f^{N-2}(x)\big) =
h\big(f^{N-1}(x)\big) = f^{N-1}(x)$ and : $g^N(x) = h\circ f^N(x) = x$.
Therefore
$x$ is $N$-periodic for $g$.
\end{enumerate}
\end{enumerate}
\hfill\carrebb

\medskip

The previous proof uses a fundamental argument : ``{\sl if $V$ is a small enough
connected open set, if
$x,y \in V$, there exists $h$ homeomorphism such that $support(h) \subset
V$ and $h(x) =
y$ }''. Let us explain why this kind of argument doesn't work for billiards~:
assume you
change a small part of the billiard table  near a point $x_0$, then you
have changed the
billiard map in a ``large'' open set, the set $U$ of all rays coming from any bounce point that is  close to $x_0$: the bounce point $x$ is close to $x_0$, but the direction of the ray may be arbitrarily chosen. The problem is that we have a fibered
system : if you
change something, you change it along a whole fiber. This problem occurs
for all
fibered problems: geodesic  flows, mechanical systems$\cdots$ and that is why we need different arguments for these cases. The problem of closing one orbit for a geodesic flow  was recently solved by L.~Rifford in \cite{Ri}.
\medskip

\noindent Without asking exactly a closing lemma, we may ask ourselves
:\medskip

``is the set of periodic points dense for a ``general'' smooth billiard table ?
''\medskip

We obtain a positive answer in the category of $C^1$-billiard tables  :\medskip

\begin{thm}\label{1}  There exists a dense $G_\delta$ subset (for the
$C^1$
topology) of the set of $C^1$-billiard tables such that, for every billiard
table  of
this $G_\delta$ subset,  the set of periodic points for the billiard map is dense in
the phase
space.

\end{thm}

Let us notice that we are unable to prove a similar result in the category of
convex $C^1$-billiard tables; this is a little surprising, because a classical
argument due to Birkhoff prove that every convex billiard table has an infinity
of periodic orbits and this result is not known for non-convex billiard tables!

We will precisely define the considered sets and topologies  in the section
2. Let us
just remark that if we perturb the billiard table in $C^1$ topology, we
perturb the
billiard map in $C^0$ topology. 

 Let us explain what is the main ingredient
of the proof
: in  \cite{B}, the authors prove that every rational polygonal billiard has a
dense set of
periodic points. Thus the main idea will be : approximating a piecewise $C^1$ curve
by a rational
polygonal one, we will create some periodic points. Using another small
perturbation, we smooth
the rational polygon at some corners and make the new periodic points non degenerate (that
means that if
we do another small perturbation, these new periodic points will still
exist). Then a
classical Baire argument is sufficient to obtain the conclusion.
\medskip

\section{\bf Definitions and results concerning the topology}\medskip

We define $\S^1 = \R/\Z$ and $\Gamma$ the set of $C^1$-maps $\gamma :
\S^1\longrightarrow
\R^2$ endowed with the $C^1$-norm (here $\| .\|_\infty$ is the norm ``sup''
associated to the usual Euclidian norm of ${\bf R}^2$)~:
$$\|\gamma\| =
\|\gamma\|_\infty +
\|\gamma'\|_\infty\; .$$ It is well-known that $(\Gamma, \|\cdot\|)$ is a
Banach space.
We denote by $d$ the associated distance.\\
We define then the set ${\cal N}$ of normal parametrisations of loops with
length $1$ :
$${\cal N} = \{\gamma \in \Gamma \;; \forall t\in {\S^1}\;,\; \|\gamma'(t)\|
= 1\}\; .$$
Then
${\cal N}$ is a closed subset of  $\Gamma$, and $({\cal N},d)$ is complete. Let us define
on ${\cal N}$ an equivalence relation $\sim$ by : $\gamma_1 \sim \gamma_2$
if  there
exists an isometry  $u : \R \longrightarrow \R$ such that $\gamma_1 =
\gamma_2 \circ
u$. If ${\cal L} = {\cal N}/\sim$, we define on ${\cal L}$ :
$$\delta(\ell_1, \ell_2) =
\inf\{d(\gamma_1, \gamma_2) ; \gamma_1 \in \ell_1, \; \gamma_2\in
\ell_2\}\; .$$ It is
easy to see that $\delta$ is a distance on ${\cal L}$ and that $({\cal L},
\delta)$ is
complete.\medskip

The subset ${\cal U}$ of ${\cal L}$ defined by : $${\cal U} = \{\ell \in
{\cal L}\;;\;
\exists \gamma \in \ell\;, \quad \forall t_1,t_2 \in {\S^1}\;,\; t_1\not= t_2
\Longrightarrow \gamma(t_1)\not= \gamma(t_2)\}$$ is an open subset of
${\cal L}$ and
therefore ${(\cal U},d)$ is a Baire space. We can identify it with the set
of $1$-dimensional compact $C^1$  submanifolds of $\R^2$ with length $1$.
\medskip

Now we can define a billiard table ; a {\it billiard table} \rm $B$ is the
closure of
the bounded connected component of the image of an element $\ell$ of ${\cal U}$ ;
then we have
$\partial B = \; Im\; \ell$ and $B$ is a simply connected $2$-dimensional
manifold with
boundary. We name
${\cal B}$ the set of billiard tables, that is too the set of  simply
connected
$2$-dimensional compact
$C^1$-submanifold with boundary of
$\R^2$ whose boundary has its length equal to 1. As the map
$\Phi : \ell \in {\cal U} \longrightarrow B\in {\cal B}$ is a bijection, we
can define a
metric $\Delta$ on ${\cal B}$ by : $$\forall(B_1,B_2) \in {\cal B}^2\;,\quad
\Delta(B_1,B_2) = \delta(\phi^{-1}(B_1)\;,\; \phi^{-1}(B_2))$$ and we know
that $({\cal
B},\Delta)$ is a Baire space ; more precisely, $({\cal B},\Delta)$ is an
open subset of a
complete metric space, and then there exists a metric $\Delta '$ on ${\cal
B}$ which is
topologically equivalent to $\Delta$ such that $({\cal B},\Delta')$ is
complete. 

Let us now give a lemma that is    interesting to understand the topology of ${\cal B}$ and that we will use later in the proof of theorem \ref{1}. \bigskip

\begin{lemma}\label{3}Let $d_H$ be the Hausdorff distance defined on
the set ${\cal
K}$ of non-empty compact subsets of $T\R^2 = \R^2\times \R^2$. For every
$B\in {\cal B}$,
let us denote  the unitary tangent fiber bundle of $\partial B$ by $K(B)$. Then
if we define on
${\cal B}$ the metric $\alpha$ by : $$\alpha(B_1,B_2) = d_H\big(K(B_1),
K(B_2)\big)$$
then $\alpha$ is topologically equivalent to $\Delta$.
\end{lemma}

\medskip

\remk
 Thus the considered topology is just the one
associated to the Hausdorff metric in the tangent fiber bundle. Let us notice
that the result is
true even if we replace the unitary tangent fiber bundle of $\partial B$ by
the unitary
tangent fiber bundle of $B$ (which is a $3$-dimensional submanifold with
boundary of
$T\R^2$)
\medskip

\medskip

In a Baire space, any countable intersection of open dense subsets is
dense. We call
{\it generic} a property which is verified by all the elements of such a
set (i.e. a
countable intersection of open dense subsets). Then any countable
intersection of
generic properties is a generic property and a generic property of a Baire
space is
satisfied by a dense subset of the  Baire space. We will
work in the
Baire space $({\cal B},\alpha)$.\medskip

Some other spaces of billiards are interesting too. The first one, named ${\cal
P}$, is the
set of piecewise $C^1$-billiard tables with length 1. Then ${\cal B}$ is a subset
of ${\cal
P}$, but there exist many elements in ${\cal P}\setminus {\cal B}$, as the
polygonal
simply connected billiards with length $1$. It is easy to define a metric
$\alpha$ on
${\cal P}$ whose restriction to ${\cal B}$ is the distance  $\alpha$ that was defined in lemma \ref{3}.   If
${\cal R}$ is
the subset of ${\cal P}$ whose elements are the rational simply connected
polygons with
length $1$ (a polygon is rational if all its angles are rational
multipliers of $\pi$),
it is easy too to see that :\medskip

\begin{lemma}\label{4}  ${\cal R}$ is dense in ${\cal P}$. Therefore,
${\cal B}$
is contained in the closure of ${\cal R}$.
\end{lemma}

We will use the  previous lemma to approximate smooth billiards  by
rational
polygonal ones, and to apply some results concerning polygonal rational
billiards.\medskip

\remk {\bf 1.} ${\cal R}$ and  ${\cal P}$  are not
Baire spaces. In these spaces, we couldn't do similar proofs. For example,
we are unable
to decide if a ``generic'' polygonal billiard has a periodic orbit.

\noindent{\bf 2.} We will denote by ${\cal P}_N$ the set   of billiards tables that  have at most $N$ corners.

 \section{Definitions and results concerning billiard maps}

Now, we consider $B\in {\cal P}$. Let $F=\{ x_1, \dots x_N\}\subset \R^2$ be the (finite) set of corners
of
$\partial B$. At every $x\in \partial B\backslash F$, we can define  the tangent
space $D(x)$ to $\partial B$ at $x$. Then we define $F(x)$ as being the set of
unitary vectors $v\in {\bf R}^2=T_x{\bf R}^2$ which are on the other side of
$D(x)$ than $B$~: in fact, $F(x)$ is a closed half-circle. Then
$\displaystyle{\Sigma=\Sigma(B)=\bigcup_{x\in \partial B\backslash F} F(x)}$ is
a 2-dimensional topological manifold with boundary. 

If $x\in \partial B\backslash F$ and $v\in D(x)\cap F(x)$, we define~:
$b(x,v)=(x,v)$.  If $x\in \partial B\backslash F$ and $v\in
F(x)\backslash D(x)$, we define~:

\begin{enumerate}

\item[$\bullet$] $w$ is the image of $v$ by the reflection of line $D(x)$;

\item[$\bullet$] $y=x+\lambda w$ where $\lambda =\inf\{ t>0; x+tw\in \partial
B\}$;

\item[$\bullet$]
$b(x, v)=(y, w)$.
\end{enumerate}
$b=b_B$ is called the {\sl billiard map}. 

\medskip

\remk Our definition of billiard map is not exactly the usual one; more
precisely, there exist two involutions $I_1$ and $I_2$ such that $b=I_1\circ
I_2$ and what is usually called the billiard map is in fact $I_2\circ
I_1=I_1\circ b\circ I_1^{-1}$; therefore the two maps are conjugated one to each
other and they have the same dynamical behavior. Moreover, our definition has
the following advantage~: $b$ is defined exactly on the set $\Sigma$.

\medskip
When $B$ is not convex, $b$ is not continuous. But it is continuous at every
point $(x, v)$ such that $p_2\circ b(x, v)$ is not tangent to $\partial B$ at
$p_1\circ b(x, v)$ (where we define~: $p_1(x, v)=x$ and $p_2(x, v)=v$). 

The following result is proved in \cite{B}~: 

\begin{thm}\label{5}(Boshernitzan, Galperin, Kruger, Troubetzkoy)~: Let $B\in
{\cal R}$ be a rational polygonal billiard; then the set of periodic points of
$b_B$ is dense in $\Sigma (B)$.
\end{thm}

We will use this result in the next section to prove    theorem \ref{1}.

\section{Proofs of  theorem \ref{1}}

Let $(U_n)_{n\in {\bf N}}$ be a countable basis of open subsets of $T^1{\bf R}^2$,
the unitary tangent fiber bundle of ${\bf R}^2$. Then we define a family
$(Q_n)_{n\in {\bf N}}$ of properties on ${\cal P}$ by~: 

\noindent{\sl ``$P\in {\cal P}$ satisfies $Q_n$ if one of the following
situations happens~:
\begin{enumerate}
\item [1.] $\Sigma (P)\cap U_n=\emptyset$;
\item [2.] $\Sigma (P)\cap U_n$ contains a periodic point for $b=b_P$''.
\end{enumerate}}

 If we prove that the set ${\cal Q}_n= \{ P\in {\cal B}; P\;{\rm verifies}\;
Q_n\}$ contains an open dense subset of ${\cal B}$, then   theorem \ref{1} is
proved.  

Thus let $n\in {\bf N}$ be fixed and ${\cal U}\subset {\cal B}$ be a
non-empty open subset of ${\cal B}$.  We have to prove that the interior of
${\cal U}\cap {\cal Q}_n$ is non-empty. Two cases are possible~:
\begin{enumerate}
\item[--] either $\forall P\in {\cal U}, \Sigma (P)\cap U_n=\emptyset$; then
${\cal U}\subset {\cal Q}_n$ and ${\cal U}\cap {\cal Q}_n ={\cal U}$ has a
non-empty interior;
\item[--] or there exists  $B_0\in {\cal U}$ such that $\Sigma (B_0)\cap
U_n\not= \emptyset$.  Then, ${\cal U}'=\{ B\in {\cal U}; \Sigma (B)\cap
U_n\not=\emptyset\}$ is a non-empty open subset of ${\cal U}$. There exists ${\cal
V}$ open subset of ${\cal P}$ such that ${\cal V}\cap {\cal B}={\cal U}'$.
Moreover, we can ask that~: $\forall P\in {\cal V}, P\cap U_n\not=\emptyset$
(because the condition ``$P\cap U_n\not=\emptyset$'' is open) and that ${\cal V}\supset\{ B\in{\cal P}; \alpha (K(B_0), K(B))<\delta\}$ for some $\delta>0$ where the distance $\alpha$ was defined in lemma \ref{3}. We have seen in
section 2 that
${\cal R}$ is dense in
${\cal P}$. Then there exists
$P_0\in {\cal R}$  such that $\alpha(K(B_0), K(P_0))<\frac{\delta }{10}$. Because $\alpha(K(B_0), K(P_0))<\frac{\delta}{10}$, $P_0$ has  $m$ corners $z_1, \dots ,z_m$ and at these corners  the distance of the two unitary tangent vectors is less that $\frac{\delta}{5}$.
As $P_0$ is a rational polygonal billiard, we
can use   theorem \ref{5}~: the billiard map $b_0$ associated to the billiard
table $P_0$ has a periodic point $(x_0, v_0)\in U_n$. 
\end{enumerate}

\noindent Because $(x_0, v_0)$ is a periodic point of $b_0$, its (finite) orbit
under
$b_0$ doesn't contain any vertex. We can smooth the billiard table near the
vertices $z_1, \dots , z_p$ without loosing the fact that $(x_0, v_0)$ is periodic; using a small
homothety (to be sure to obtain a length equal to 1), we obtain a $C^2$-billiard
table $P_1$  such that $(x_0, v_0)\in
U_n\cap\partial P_1$ is periodic for the billiard map $b_1$ associated to $P_1$. Because  $\alpha(K(B_0), K(P_0))<\frac{\delta}{10}$ and because   the distance of the two unitary tangent vectors at these corners of $P_0$ is less that $\frac{\delta}{5}$, we can ask that $\alpha(K(B_0), K(P_1))<\delta$ and then that $P_1\in {\cal V}\cap {\cal B}={\cal U}'$.

We have not finished the proof of   theorem \ref{1} because we have find $P_1\in
{\cal U}\cap {\cal Q}_n$ but we don't know if ${\cal U}\cap {\cal Q}_n$ has a
non-empty interior. The next idea is then to perturb $P_1$ in such a way that
$(x_0, v_0)$ becomes stably periodic.

To do that, we recall some results contained in \cite{KT}; we will call $(x_1,
v_1)=b_1(x_0, v_0)$, \dots, $(x_\tau , b_\tau )=b_1^\tau (x_0, v_0)=(x_0, v_0)$
the points of the (periodic) orbit of $(x_0, v_0)$ under $b_1$.  Let us name
$(y_i, w_i)$ coordinates near $(x_i, v_i)$. In a neighbourhood of $(x_i, v_i)$,
the map $((y_i, w_i)\rightarrow (y_i, z_i))$ where $z_i=p_1\circ b_1(y_i, w_i)$
is a $C^1$-diffeomorphism. Then, $(x_0, \dots, x_{\tau -1})$ (abbreviation  for
$((x_0, x_1), (x_1, x_2), \dots , (x_{\tau -1}, x_0))$ in these coordinates) is a
periodic orbit if and only if it is  a critical point of the
$C^2$ length function defined by (we note $y_\tau =y_0$ in this case)~:
$$\ell (y_0, \dots , y_{\tau -1})=\sum_{i=1}^\tau \| y_i-y_{i-1}\| .$$
At such a point $(x_0, \dots , x_{\tau -1})$, we define (eventually in charts)~:
$$a_i=\frac{\partial^2\ell}{\partial y_i^2}(x_0, \dots , x_{\tau -1})\;{\rm
and}\; b_i=\frac{\partial^2\ell}{\partial y_i\partial y_{i+1}}(x_0, \dots ,
x_{\tau -1}).$$
Then the Hessian of $\ell$ is~:
\begin{enumerate}
\item[$\bullet$] if $\tau =2$~: $H_2=\begin{pmatrix}a_1&b_1+b_2\\ b_2+b_1 &
a_2\end{pmatrix}$;
\item[$\bullet$] il $\tau >2$~:  $H_\tau =\begin{pmatrix}a_1&b_1&0&\dots&b_\tau\\
b_1 & a_2&b_2&\dots &0\\0&b_2&a_3&\dots&0\\ \dots&&\dots&&\dots \\
b_\tau&0&0&\dots&a_\tau\end{pmatrix}$.
\end{enumerate}

An easy  calculus made in \cite{KT} shows that~:
\begin{enumerate}
\item[--] $b_i$ depends only on $x_i$, $x_{i+1}$ and the tangent space to
$\partial P_1$ at $x_i$ and $x_{i+1}$;
\item[--] if the tangent space to $\partial P_1$ at $x_i$ is fixed, $a_i$
depends linearly on the curvature of $\partial P_1$ at $x_i$, and this
dependance is effective ($a_i$ is not constant).
\end{enumerate}

As it is easy to change the curvature of a curve near a point without changing
the tangent space at this point or far away this point, we can perturb $P_1$ in
$P_2\in {\cal B}\cap {\cal U}$ that is $C^2$ and such that~: $\det H_\tau \not=
0$; indeed, $\det H_\tau$ is a non constant polynomial function in $a_1$, \dots ,
$a_\tau$. For this new billiard $P_2$, $(x_1, \dots , x_{\tau -1})$ is a
non-degenerate critical point of the $C^2$ function $\ell_2$. We have~:
\begin{lemma}\label 6
There exists a neighbourhood $W$ of $\ell_2$ in the $C^1$ topology such that
every element of $W$ has a critical point near $(x_0, \dots , x_{\tau -1})$.
\end{lemma}

\remk the previous lemma is more known in the case of the $C^2$-topology. In the
case of the $C^2$-topology, we can even ask that the critical point near $(x_0
\dots , x_{\tau -1})$ is unique. In the $C^1$-topology, we may have an infinity
of such critical points, but we know that at least one of these critical points
exists. Let us explain why.\\
 If we consider $\ell$ that is $C^2$ but just $C^1$ close to $\ell_2$:
\begin{enumerate}
\item[--] we use the existence of an isolating block $B$ for ${\rm grad}\ell$, which
is stable by $C^1$-perturbation ( see
\cite{C}); this implies the existence of one positive (or negative for a minimum) 
orbit for the flow $(\varphi_t)$ of $\ell$ which stays in
$B$;
\item[--] but if $\ell$ has no critical point in $B$, this is impossible, because
there exists a constant $k>0$ such that~:
$\forall x\in B, {d\ell\circ \varphi_t\over dt}(x)=\| {\rm grad}\ell(x)\|^2\geq
k$, and then $\displaystyle{\lim_{t\rightarrow +\infty} \ell\circ
\varphi_t(x)=+\infty}$ (and $\displaystyle{\lim_{t\rightarrow -\infty} \ell\circ
\varphi_t(x)=-\infty}$) and thus the orbit leaves $B$.

\end{enumerate}
If $\ell$ is not $C^2$, we cannot use the same argument because the gradient flow is not defined. Let us assume that $\ell$ is $C^1$ close to $\ell_2$ and has no critical point in the isolating block. By using a convolution, we can approximate $\ell$ in $C^1$ topology by a smooth $\ell_1$ that has no critical point in the isolating block, and this contradicts what we explained for $C^2$ functions $\ell$. 

Now, we can finish the proof of  theorem \ref{1}~: every billiard table $B$
$C^1$-close to $P_2$ has its length function $\ell$ which is $C^1$-close to
$\ell_2$; thus the associated billiard map has at least one periodic orbit close to
$((x_0, v_0), \dots , (x_{\tau -1}, v_{\tau -1}))$ and therefore a periodic point
in $U_n$. Then a whole neighbourhood of $P_2$ in ${\cal B}$ is in ${\cal Q}_n$ and
then the interior of ${\cal U}\cap {\cal Q}_n$ is non-empty. \hfill\carrebb 

\medskip
\remk  
\noindent{\bf 1.} The problem of convex billiards is less easy. Of course, we can
prove that their set is a Baire set and we can approximate any convex $C^1$-billiard
table by a convex rational polygonal one having a periodic orbit; we can smooth
this billiard to obtain a convex one;  but if the periodic orbit of the rational
polygonal billiard has more than one bounce point on every side, you cannot perturb
it in such a way that you obtain a strictly convex
$C^2$-billiard table having the same periodic orbit. And if the smooth, the  billiard
table that we obtain is convex but not strictly convex, we cannot ask, when we
perturb it in the
$C^2$-topology to obtain a stable periodic orbit, that the new billiard table is
convex.

\noindent {\bf 2.} We know that a generic element of ${\cal B}$ is topologically
transitive (see \cite{G} for a proof in a slightly different topology, but easily
transposable in our topology by using \cite{Ka}).  A dynamical system having a dense subset of
periodic points and being topologically transitive is called chaotic (see
\cite{BB}); we have then prove that the billiard map associated to a generic $C^1$
billiard table is chaotic.
\bigskip

\bibliography{test}

\bibliographystyle{alpha}

\noindent {\bf{Bibliography.}}


\end{document}